\documentclass[a4paper,12pt]{article}
\usepackage[T1]{fontenc}
\usepackage[utf8]{inputenc}
\usepackage{graphicx}
\usepackage{fancyhdr}
\usepackage{float}
\usepackage{xcolor}
\usepackage{authblk}
\usepackage{mathrsfs}
\usepackage{empheq}
\usepackage[hyphens]{url}
\usepackage{hyperref} 
\usepackage[]{breakurl}
\usepackage[]{amsthm}
\usepackage{tikz}
\usepackage{bm}
\usepackage{amsmath}
\usepackage{amssymb}
\usepackage{url}
\usepackage{cite}
\usepackage{multirow,booktabs}

\usepackage{geometry}
\begin{document}

\title{Estimates of the minimum of the Gamma function using the Lagrange inversion theorem and the Fa\`a di Bruno formula}

\author[$\dagger$]{Jean-Christophe {\sc Pain}$^{1,2,}$\footnote{jean-christophe.pain@cea.fr}\\
\small
$^1$CEA, DAM, DIF, F-91297 Arpajon, France\\
$^2$Universit\'e Paris-Saclay, CEA, Laboratoire Mati\`ere en Conditions Extr\^emes,\\ 
F-91680 Bruy\`eres-le-Ch\^atel, France
}

\maketitle

\begin{abstract}
In this article we derive, using the Lagrange inversion theorem and applying twice the Fa\`a di Bruno formula, an expression of the minimum of the Gamma function $\Gamma$ as an expansion in powers of the Euler-Mascheroni constant $\gamma$. The result can be expressed in terms of values the Riemann zeta function $\zeta$ of integer arguments, since the multiple derivative of the digamma function $\psi$ evaluated in $1$ is precisely proportional to the zeta function. The first terms (up to $\gamma^6$) were provided in order to address the convergence of the series. Applying the Lagrange inversion theorem at the value $3/2$ yields more accurate results, although less elegant formulas, in particular because the digamma function evaluated in $3/2$ does not simplify.
\end{abstract}

\section{Introduction}

Let us assume that $z$ is defined as a function of $w$ by an equation of the form $w=f(z)$, where $f$ is analytic at a point a and $f'(a)\ne 0$. Then, according the the Lagrange theorem \cite{Lagrange1770,Whittaker1990} it is possible to invert or solve the equation for $z$, expressing it in the form $z=h(w)$ given by a power series
\begin{equation}
h(w)=a+\sum _{n=1}^{\infty}h_{n}{\frac {(z-f(a))^{n}}{n!}},
\end{equation}
where
\begin{equation}
h_{n}=\lim _{z\to a}{\frac {d^{n-1}}{dz^{n-1}}}\left[\left({\frac {z-a}{f(z)-f(a)}}\right)^{n}\right].
\end{equation}
Let us set 
\begin{equation}
f(z)=\psi(z)=\frac{\Gamma'(z)}{\Gamma(z)},
\end{equation}
where $\psi(z)$ represents the digamma function. In order to find the minimum of the Gamma function $\Gamma(z)$, and using the well-known properties of the latter \cite{Artin1964}, we can apply the above mentioned Lagrange theorem to $\psi(z)$ \cite{Campbell1966}. The latter can be written, for $n\ne -1, -2, -3, \cdots$, as
\begin{equation}
    \psi(z)=-\gamma+\sum_{n=0}^{\infty}\left(\frac{1}{n+1}-\frac{1}{n+z}\right),
\end{equation}
where $\gamma$ is the usual Euler-Mascheroni constant, and has clearly no poles for $z>0$. Setting
\begin{equation}
    \mathcal{L}(z)=\frac{z-a}{\psi(z)-a},
\end{equation}
we get that the minimum of $\Gamma$ can be expanded as
\begin{equation}\label{mini}
    z_m=\left.a+\sum_{n=1}^{\infty}\frac{\left[-\psi(a)\right]^n}{n!}\,\frac{d^{n-1}}{dz^{n-1}}\left[\mathcal{L}(z)\right]^n\right|_{z=a}.
\end{equation}

In section \ref{sec2}, we provide a general formula for the minimum of the Gamma function as a series expansion, obtained by combining Eq. \ref{mini} with the Fa\`a di Bruno formula. The first terms up to $\gamma^6$, in the $a=1$ case, are explicitly given in section \ref{sec3}, as well as an expansion obtained at $a=3/2$.

\section{General expression of the minimum as an infinite series}\label{sec2}

The Fa\`a di Bruno formula \cite{Arbogast1800,Bruno1855,Bruno1857,Comtet1974,Johnson2002,Stanley1999,Andrews1976,Yang2000} gives
\begin{equation}
\frac{d^{n-1}}{dz^{n-1}}\left[\mathcal{L}(z)\right]^n=(\mathcal{L}^n)^{(n-1)}(z)=\sum_{\ell_1+\ell_2+\cdots+\ell_{n}=n-1}\binom{n-1}{\ell_1,\ell_2,\cdots,\ell_n}\,\prod_{i=1}^{n}\mathcal{L}^{(\ell_i)}(z),
\end{equation}
where
\begin{equation}
\binom{n-1}{\ell_1,\cdots,\ell_n}=\frac{(n-1)!}{\ell_1!\,\ell_2!\,\cdots\,\ell_n!}
\end{equation}
is the usual multinomial coefficient. We have also, thanks to the Leibniz rule for the multiple derivative of a product:
\begin{equation}
    \mathcal{L}^{(\ell_i)}(a)=\left.\frac{\partial^{\ell_i}}{\partial z^{\ell_i}}\left[\frac{z-a}{\psi(z)-a}\right]\right|_{z=a}=\left.\ell_i\,\frac{\partial^{\ell_i-1}}{\partial z^{\ell_i-1}}\left[\frac{1}{\psi(z)-a}\right]\right|_{z=a}.
\end{equation}

This enables us to write, according to the Fa\`a di Bruno formula again:
\begin{equation}
    \left.\frac{d^{\ell_i-1}}{dz^{\ell_i-1}}\left[\frac{1}{\psi(z)-a}\right]\right|_{z=a}=\left.\frac{(\ell_i-1)!}{\left[\psi(z)-a\right]^{\ell_i}}\sum {\frac {(-1)^{\ell_i-1-k_{0}}(\ell_i-1-k_{0})!}{\displaystyle \prod _{j=1}^{\ell_i-1}(j!)^{k_{j}}\;k_{j}!}}\prod _{j=0}^{\ell_i-1}\left[\psi^{(j)}(z)-a\,\delta_{j,0}\right]^{k_{j}}\right|_{z=a},
\end{equation}
where $k_{0}=k_{2}+2k_{3}+\dots +(\ell_i-2)k_{\ell_i-1}$ is such that $n=k_{0}+k_{1}+\dots +k_{\ell_i-1}=k_{1}+2k_{2}+\dots +(\ell_i-1)k_{\ell_i-1}$. $\delta_{m,n}$ is the usual Kronecker symbol.

Equation \ref{mini} can be put in the form
\begin{align}\label{mini2}
    z_m=&a+\sum_{n=1}^{\infty}\frac{\left[-\psi(a)\right]^n}{n!}\,\sum_{\ell_1+\ell_2+\cdots+\ell_{n}=n-1}\binom{n-1}{\ell_1,\ell_2,\cdots,\ell_n}\,\prod_{i=1}^{n}\ell_i\nonumber\\
    &\left.\times\frac{(\ell_i-1)!}{\left[\psi(z)-a\right]^{\ell_i}}\sum\frac{(-1)^{\ell_i-1-k_{0}}(\ell_i-1-k_{0})!}{\prod_{j=1}^{\ell_i-1}(j!)^{k_{j}}\;k_{j}!}\prod_{j=0}^{\ell_i-1}\left[\psi^{(j)}(a)-a\,\delta_{j,0}\right]^{k_{j}}\right|_{z=a},
\end{align}
which is the main result of the present work. We have actually, for any positive integer $p$:
\begin{equation}\label{rel}
    \psi^{(n)}(s)=(-1)^{n+1}\,n!\,\zeta(n+1,s)
\end{equation}
where
\begin{equation}
    \zeta(u,v)=\sum_{k=0}^{\infty}\frac{1}{(k+v)^u}
\end{equation}
is the Hurwitz zeta function and $\zeta(u,1)=\zeta(u)$ the usual zeta function.

In the general case (i.e., for any value of $a$ in the interval considered here), Eq. (\ref{mini2}) becomes therefore
\begin{align}
    z_m=&a+\sum_{n=1}^{\infty}\frac{\left[-\psi(a)\right]^n}{n!}\,\sum_{\ell_1+\ell_2+\cdots+\ell_{n}=n-1}\binom{n-1}{\ell_1,\ell_2,\cdots,\ell_n}\,\prod_{i=1}^{n}\ell_i\nonumber\\
    &\left.\times\frac{(\ell_i-1)!}{\left[\psi(z)-a\right]^{\ell_i}}\sum\frac{(-1)^{\ell_i-1-k_{0}}(\ell_i-1-k_{0})!}{\prod_{j=1}^{\ell_i-1}(j!)^{k_{j}}\;k_{j}!}\prod_{j=0}^{\ell_i-1}\left[(-1)^{j+1}\,j!\,\zeta(j+1,a)-a\,\delta_{j,0}\right]^{k_{j}}\right|_{z=a},
\end{align}
with $k_{0}=k_{2}+2k_{3}+\dots +(\ell_i-2)k_{\ell_i-1}$ is such that $n=k_{0}+k_{1}+\dots +k_{\ell_i-1}=k_{1}+2k_{2}+\dots +(\ell_i-1)k_{\ell_i-1}$.

\section{The first terms and an insight into the convergence}\label{sec3}

\subsection{Case $a=1$: series expansion in powers of $\gamma$ (up to $\gamma^6$)}\label{subsec31}

In the case where $a=1$, Eq. (\ref{rel}) becomes
\begin{equation}
    \psi^{(n)}(1)=(-1)^{n+1}\,n!\,\zeta(n+1)
\end{equation}
and one has
\begin{equation}
    z_m=\left.1+\sum_{n=1}^{\infty}\frac{\left[-\psi(1)\right]^n}{n!}\,\frac{d^{n-1}}{dz^{n-1}}\left[\mathcal{L}(z)\right]^n\right|_{z=1},
\end{equation}
or
\begin{align}
    z_m=&1+\sum_{n=1}^{\infty}\frac{\left[-\psi(1)\right]^n}{n!}\,\sum_{\ell_1+\ell_2+\cdots+\ell_{n}=n-1}\binom{n-1}{\ell_1,\ell_2,\cdots,\ell_n}\,\prod_{i=1}^{n}\ell_i\nonumber\\
    &\left.\times\frac{(\ell_i-1)!}{\left[\psi(z)-1\right]^{\ell_i}}\sum\frac{(-1)^{\ell_i-1-k_{0}}(\ell_i-1-k_{0})!}{\prod_{j=1}^{\ell_i-1}(j!)^{k_{j}}\;k_{j}!}\prod_{j=0}^{\ell_i-1}\left[(-1)^{j+1}\,j!\,\zeta(j+1)-a\,\delta_{j,0}\right]^{k_{j}}\right|_{z=1}
\end{align}
with $k_{0}=k_{2}+2k_{3}+\dots +(\ell_i-2)k_{\ell_i-1}$ is such that $n=k_{0}+k_{1}+\dots +k_{\ell_i-1}=k_{1}+2k_{2}+\dots +(\ell_i-1)k_{\ell_i-1}$.

We obtain, using the above expression and the computer algebra system Mathematica \cite{Mathematica}, the following expression :
\begin{equation}
    z_m=1+r_1\,\gamma+r_2\,\gamma^2+r_3\,\gamma^3+r_4\,\gamma^4+r_5\,\gamma^5+r_6\,\gamma^6+\cdots,
\end{equation}
where
\begin{equation}
    r_1=\frac{1}{\zeta^2(2)},
\end{equation}
\begin{equation}
    r_2=\frac{\zeta(3)}{\zeta^3(2)},
\end{equation}
\begin{equation}
    r_3=\left[\frac{2\zeta^2(3)}{\zeta(2)}-\zeta(4)\right]\frac{1}{\zeta^4(2)},
\end{equation}
\begin{equation}
    r_4=\left[-2\zeta(2)\zeta(3)+2\frac{\zeta^2(3)}{\zeta(4)}+\zeta(5)\right]\frac{1}{\zeta^5(2)},
\end{equation}
\begin{equation}
    r_5=\left[-\frac{42}{5}\zeta^2(3)+\frac{16}{5}\frac{\zeta^4(3)}{\zeta(6)}+6\frac{\zeta(3)\zeta(5)}{\zeta(2)}+\frac{11}{10}\zeta(6)\right]\frac{1}{\zeta^6(2)},
\end{equation}
and
\begin{equation}
    r_6=\left[\frac{36}{5}\zeta(3)\zeta(4)-168\frac{\zeta^2(3)}{\zeta(2)}+\frac{144}{25}\frac{\zeta^5(3)}{\zeta(8)}-\frac{14}{5}\zeta(2)\zeta(5)+\frac{56}{5}\frac{\zeta^2(3)\zeta(5)}{\zeta(4)}+\zeta(7)\right]\frac{1}{\zeta^7(2)}.
\end{equation}
The values obtained for expansions up to sixth order are given in table \ref{tab1}. We can see that the convergence is rather slow, and characteristic of an asymptotic expansion.

\begin{table}
\centering
    \begin{tabular}{lc}\hline\hline
    Truncated series & Value of the root\\\hline
    \addlinespace
    1+$r_1\,\gamma$ & 1.213324688\\
    \addlinespace
    1+$r_1\,\gamma$+$r_2\,\gamma^2$ & 1.303306712\\
    \addlinespace
    1+$r_1\,\gamma$+$r_2\,\gamma^2$+$r_3\,\gamma^3$ & 1.433026242\\
    \addlinespace
    1+$r_1\,\gamma$+$r_2\,\gamma^2$+$r_3\,\gamma^3$+$r_4\,\gamma^4$ & 1.465429144\\
    \addlinespace
    1+$r_1\,\gamma$+$r_2\,\gamma^2$+$r_3\,\gamma^3$+$r_4\,\gamma^4$+$r_5\,\gamma^5$ & 1.471535623\\
    \addlinespace
    1+$r_1\,\gamma$+$r_2\,\gamma^2$+$r_3\,\gamma^3$+$r_4\,\gamma^4$+$r_5\,\gamma^5$+$r_6\,\gamma^6$ & 1.472388063 \\\addlinespace\hline\hline
    \end{tabular}
\caption{Results obtained when the Lagrange inversion is performed at $a=1$. the exact solution is 1.4616321449683623413...}\label{tab1}
\end{table}

\subsection{Case $a=3/2$: a better choice, but more complex expressions}

In the case where $a=3/2$, we obtain
\begin{align}
z_m=&\frac32+\frac{2}{\pi^2-8}(-2+\gamma+\ln4)+\frac{8(7\zeta(3)-8)}{(\pi^2-8)^3}(-2+\gamma+\ln4)^2\nonumber\\
&+\left[\frac{64(7\zeta(3)-8)^3}{(\pi^2-8)^5}-\frac{8(\pi^2-96)}{3(\pi^2-8)^4}\right](-2+\gamma+\ln4)^3+\cdots
\end{align}
One gets, setting
\begin{equation}
    q_1=\frac{2}{\pi^2-8}(-2+\gamma+\ln4),
\end{equation}
\begin{equation}
    q_2=\frac{2}{\pi^2-8}(-2+\gamma+\ln4)+\frac{8(7\zeta(3)-8)}{(\pi^2-8)^3}(-2+\gamma+\ln4)^2,
\end{equation}
and
\begin{equation}
    q_3=\left(\frac{64(7\zeta(3)-8)^3}{(\pi^2-8)^5}-\frac{8(\pi^2-96)}{3(\pi^2-8)^4}\right)\,\eta^3,
\end{equation}
where $\eta=-2+\gamma+\ln4$, the results displayed in table \ref{tab2} for the three successive approximations 1+$q_1$, 1+$q_2$ and 1+$q_3$.

\begin{table}
\centering
    \begin{tabular}{lc}\hline\hline
    Truncated series & Value of the root\\\hline
   \addlinespace$\frac{3}{2}$+$q_1$ & 1.460965032\\ 
   \addlinespace$\frac{3}{2}$+$q_1$+$q_2$ & 1.461640502\\
   \addlinespace$\frac{3}{2}$+$q_1$+$q_2$+$q_3$ & 1.461632068\\\addlinespace\hline\hline
   \end{tabular}
\caption{Results obtained when the Lagrange inversion is performed at $a=3/2$.}\label{tab2}
\end{table}

The convergence is much more faster than for the case $a=1$ \cite{Stack}, but in the $a=3/2$ case unfortunately, $\zeta(u,3/2)$ does not yield, to our knowledge, to a known special function...

\section{Conclusion}

We obtained, using the Lagrange inversion theorem (at the value $1$) and applying twice the Fa\`a di Bruno formula, an expression of the minimum of the Gamma function as an expansion in powers of the Euler-Mascheroni constant. The expression involves only values of the Riemann zeta function $\zeta$ of integer arguments, due to the fact that the multiple derivative of the digamma function evaluated in $1$ is actually proportional to the zeta function. The first terms (up to $\gamma^6$) were provided in order to address the convergence of the series. Applying the Lagrange inversion theorem at the value $3/2$ gives much more accurate results, although less elegant formulas, since the multiple derivative of the digamma function evaluated in $3/2$ does not reduce, to our knowledge, to a usual special function.

\end{document}